\documentclass[12pt,a4paper]{amsart}
  \usepackage{amsmath,amssymb,amsthm,amscd,latexsym,graphicx,tikz}
 \usepackage{epigraph}
 \usepackage[OT2,OT1]{fontenc}
  \usepackage{enumerate}

\usepackage{graphicx}
\usepackage[all]{xy}

  \entrymodifiers={!!<0pt,0.7ex>+}

 \renewcommand{\U}{\mathbf{U} }

\newcommand{\carac}{\mathrm{char}}

 \newcommand{\obs}{\mathrm{Obs}}

 \renewcommand{\to}{\longrightarrow}

 \newcommand{\C}{\mathbb{C}}
 \newcommand{\F}{\mathbb{F}}

 \newcommand{\Q}{\mathbb{Q}}
 
 \newcommand{\Z}{\mathbb{Z}}

 \newcommand{\GL}{\mathbf{GL}}

 \newcommand{\End}{\mathbf{End}}

 \newcommand{\Br}{\mathrm{Br}}

\newcommand{\Gal}{\mathrm{Gal}}

\renewcommand{\ni}{\noindent}

 \newcommand{\B}{\mathbf{B}}

 \newcommand{\Lie}{\mathrm{Lie}}

\renewcommand{\Im}{\mathrm{Im}}
 
\newcommand{\nat}{\mathrm{nat}}
\newcommand{\Ker}{\mathrm{Ker}}

\newcommand{\W}{\mathbf{W}}

\newcommand{\frob}{\mathrm{frob}}

\newcommand{\vers}{\mathrm{vers}}

 \theoremstyle{plain}
 \newtheorem{thm}{Theorem}[section]
   \newtheorem*{thm*}{Theorem}
 \newtheorem{defi}[thm]{Definition}

 \newtheorem{prop}[thm]{Proposition}

 \newtheorem{lem}[thm]{Lemma}

 \newtheorem*{theorem-non}{Bloch-Kato conjecture, an equivalent formulation}
\newtheorem*{thmA*}{Theorem A}
\newtheorem*{thmB*}{Theorem B}
\newtheorem*{thmC*}{Theorem C}
\newtheorem*{thmD*}{Theorem D}

 \theoremstyle{remark}
 \newtheorem{rem}[thm]{Remark}

 \newtheorem{qu}[thm]{Question}

 \newtheorem{ex}[thm]{Example}

 \newenvironment{dem}{{\bf Proof.}}{\hfill$\square$}
 
\begin{document}
\title{Lifting Galois representations via Galois cohomology}

\author{Mathieu Florence}
\address{Université Paris Cité and Sorbonne Université, CNRS, IMJ-PRG, F-75005 Paris, France }
\email{mathieu.florence@imj-prg.fr}





\begin{abstract}
    
This is a survey on some recent  developments about lifting Galois representations.
\end{abstract}
\maketitle
\tableofcontents
\newpage
\section{Introduction.}
\ni This paper surveys  recent progress on lifting  Galois representations. We focus on approaches that solely rely on  Galois cohomology: \cite{CDF}, \cite{DCF0}, \cite{KL}, \cite{MeSc1} and \cite{MeSc2}. For an overview of this widely investigated topic, also beyond the scope of this survey, see for instance \cite{CDF} or  \cite{MeSc1}. Observe that sections \ref{Heisenbergpodd} and \ref{secnonKum} actually contain new material. In section 6, emphasis is laid on results of \cite{CDF}, on which I gave a talk at the conference held in Ottawa in June 2024, to celebrate J\'an Min\'a\v{c}'s 71th birthday. \\Dear J\'an, it is a pleasure to know you. I wish you lots of enjoyable mathematics in your life!
\section{General facts on lifting problems.}
\ni Let $p$ be a prime. In this text, $G$ always denotes a profinite group.

\begin{defi}
     A homomorphism from $G$ to an abstract group is said to be continuous, if its kernel is open. 
\end{defi}
\begin{defi}
 Let $R$ be a commutative ring.  An $(R,G)$-module is a finite locally free $R$-module, equipped with an $R$-linear $G$-action, whose kernel is open in $G$.    
\end{defi}
\ni Consider the following  general question.

\begin{qu}\label{QuGen}
   Let $G$ be a profinite group. Let $\Gamma$ be an affine and smooth $\Z$-group scheme.  Consider a  continuous representation  \[\rho: G \to \Gamma(\F_p).\]   Does $\rho$ admit a continuous lift  \[\rho_2: G \to \Gamma(\Z/p^2)? \]
\end{qu}

  \ni When $G$ is an absolute Galois group  and $\Gamma:=\GL_d$, Question \ref{QuGen} yields the main problem of interest in this paper- ubiquitous in arithmetic.  

\begin{qu}\label{qugallift}
   Let $F$ be a field, with a separable closure $F_s$. \\For some integer $d \geq 1$, consider  a  (continuous) Galois representation \[\rho: \Gal(F_s/F) \to \GL_d(\F_p).\]  Does $\rho$ admit a lift, to some  \[\rho_2: \Gal(F_s/F)  \to \GL_d(\Z/p^2)? \]
\end{qu}

\subsection{Variants and generalisations.\\}
\ni  Denote by $\B_d \subset \GL_d$ the Borel subgroup of upper triangular matrices, and by $\U_d \subset \B_d$ its unipotent radical, consisting of strictly upper triangular matrices. There is a `triangular' (resp. 'strictly triangular') variant of Question \ref{qugallift}, replacing $\GL_d$ by $\B_d$ (resp. by $\U_d$). In some cases, these triangular variants are easier to handle, by recursive lifting algorithms. This is the main guideline of \cite{CDF}.\\

In the formulation of Question \ref{QuGen}, the field of coefficients of  $\rho$ is $\F_p$. One may replace it by an arbitrary field $k$
 of characteristic $p$, and accordingly, replace $\Z/p^2$ by $\W_2(k)$ (truncated Witt vectors of length two). Depending on the results aimed at, this may, or may not, make a difference. For instance, the proof of Proposition \ref{liftp2d3} (liftability of two-dimensional Galois representations) is insensitive to $k$. On the other hand, for $p=2$,  both  statement and proof of Theorem \ref{thmnonlift} (that includes the non-liftability of the generic Galois representation of dimension $5$ over $k=\F_2$) highly depend on $k$.\\
 
One may also search for liftings  modulo $p^r$, for $r \geq 3$. This topic is addressed in Section \ref{secspeclift}, for absolute Galois groups of local fields (and more generally, for the so-called $p$-manageable profinite groups).

\subsection{Cohomological obstructions to lifting.\\}\label{secobs}
Place ourselves in the context of Question \ref{QuGen}. Endow the Lie algebra $\Lie(\Gamma)$ with its natural (adjoint)  linear $\Gamma$-action. Denote by \[\End_{\F_p}(\rho) := \rho^*(\Lie_{\F_p}(\Gamma))\] the Lie algebra of $\Gamma$ over $\F_p$, considered as a representation of $G$, by group-change via $\rho.$
    By the general formalism of (group) cohomology, the obstruction to the existence of $\rho_2$, is a natural class \[\obs(\rho_2) \in H^2(G,\End_{\F_p}(\rho)).\] This class can be described explicity via group extensions; see e.g. section 2.2. of \cite{MeSc1} for the cases $\Gamma=\GL_d,\B_d$.

\begin{rem}(Arbitrary field of coefficients.)\\
Let $k$ be a field of characteristic $p$.  Denote by \[\frob: k \xrightarrow{x \mapsto x^p} k\] the frobenius endomorphism of $k$. For a $k$-vector space $V$, denote by \[V^{(1)}:=V \otimes_{\frob} k\] its frobenius twist. Consider a continuous representation \[\rho: G \to \Gamma(k).\]  Then, the  obstruction to lifing $\rho$ to \[\rho_2: G \to \Gamma(\W_2(k)), \] is  a natural class \[\obs(\rho_2) \in H^2(G,\End_k(\rho)^{(1)}).\]
The  frobenius twist, invisible when $k=\F_p$, is actually essential.
\end{rem}

\subsection{Subgroups of prime-to-$p$ index, reduction to the triangular case.\\}\label{remrescor}
Let $G$ be a profinite group, and consider a representation   \[\rho: G \to \GL_d(\F_p). \] There is an open subgroup $G_0 \subset G$, of prime-to-$p$ index, such that the restriction $\rho_{\vert G_0}$ is upper triangular, up to conjugation. One may thus assume that   $\rho_{\vert G_0}$ reads as\[\rho_{\vert G_0}: G_0 \to  \B_d(\F_p) \subset \GL_d(\F_p). \] 
 Assume that $ \rho_{\vert G_0}$ lifts, to a representation   \[ G_0 \to \B_d(\Z/p^2). \]  Then $\rho$ lifts, to a representation  \[\rho_2: G \to \GL_d(\Z/p^2). \] The same result holds, replacing $\B_d \subset \GL_d$, by $\U_d \subset \B_d$.\\ The proof is  a classical restriction/corestriction argument, to a pro-$p$-Sylow $G_p \subset G$. See for instance \cite{MeSc2}, Lemma 2.6.

\begin{rem}
This restriction/corestriction argument  typically does not apply to lifting modulo higher powers of $p$. Precisely, if one replaces $\Z/p^2$ by $\Z/p^3$, it is likely that the analoguous result is false- though I do not have a counter-example in mind. In general, there is the following fact, that one should relate to Hensel's Lemma. For some $n\geq 1$, consider a representation $\rho_n: G \to \Gamma(\Z/p^n)$. Then, liftability of $\rho_n$, to a representation $\rho_{2n}: G \to \Gamma(\Z/p^{2n})$, is an abelian problem: it is obstructed by a natural class \[\obs(\rho_{2n}) \in H^2(G,\End_{\Z/p^n}(\rho)).\] However, lifting $\rho_n$, to $\rho_{2n+1}: G \to \Gamma(\Z/p^{2n+1})$ is not an abelian question:  there is no natural cohomology class (with values in a $G$-module) that obstructs it.
\end{rem}

\subsection{A simple descent Lemma.}
 \begin{lem}\label{lemtech}
     Let $G$ be a (profinite) group, and let  $V,V'$ be  $(k,G)$-modules. Let $l/k$ be an extension of fields of characteristic $p$. The following holds. \begin{enumerate}
         \item{If $\left( V \oplus V' \right)$ lifts to a  $(\W_2(k),G)$-module, then so does $V$.} \item{If the  $(l,G)$-module $V \otimes_k l$ lifts to a  $(\W_2(l),G)$-module, then $V$ lifts to a $(\W_2(k),G)$-module $V_2$.}
         
     \end{enumerate}
 \end{lem}
 \begin{dem}
     Item (1) is \cite{DCF0}, Lemma 3.4. Item (2) is classical. It follows from the facts, that the obstruction to the existence of $V_2$, reading as \[\obs(V_2) \in H^2(G,\End_k(V)^{(1)})\] (see section \ref{secobs}) is compatible to base-change,  and that the injection\[\End_k(V)^{(1)} \hookrightarrow  \End_l(V\otimes_k l)^{(1)}= \End_k(V^{(1)})  \otimes_k l\] has a $G$-equivariant retraction- provided by the choice of a $k$-linear retraction of the inclusion $k \hookrightarrow l$.
 \end{dem}

\section{A non-liftable Heisenberg-Galois representation.}\label{Heisenbergpodd}

\ni In this survey, we mostly (but not exclusively) focus on lifting mod $p$ Galois representations. By definition, they are representations of $G:=\Gal(F_s/F)$, where $F$ is a field, with values in $\GL_d(k)$, where $k$ is a field of characteristic $p$. Such a representation that takes values in  $\B_d(k)$ resp. $\U_d(k)$, is called triangular, resp. strictly triangular.

\ni We begin with a simple significant example, taken from an unpublished earlier version of \cite{F2}. It was removed from  recent versions.\\

\ni For $p \geq 3$, we give an elementary example of a field $F$, containing $\C$, such that the following natural arrow  is not surjective: \[H^1(\Gal(F_s/F), \mathbf U_3(\Z/p^2)) \to H^1(\Gal(F_s/F), \mathbf U_3(\F_p)).\]
Thus,  all  mod $p$ ``Heisenberg-Galois'' representations do not lift mod $p^2$. \\

\ni Start with a field $F$, containing $\C$. Set $G$ to be its absolute Galois group. For each $n \geq 1$, use $e^{\frac{2\pi i} n } \in F$ to identify $\mu_n$ to $\Z/n$, as finite $G$-modules. Pick $x,y \in F^\times$.\\

\ni By Kummer theory, there are two classes $$(x)_p,(y)_p \in H^1(F,\mu_p), $$ respectively associated to extensions of $(\F_p,G)$-modules \[\mathcal E_x: 0 \to  \F_p =\mu_p \to E_x \to \F_p \to 0\]  and \[\mathcal E_y: 0 \to  \F_p=\mu_p \to E_y \to \F_p \to 0.\] These give rise to group homomorphisms \[\rho_x: G \to \mathbf U_2(\F_p) =\F_p\] and \[\rho_y: G \to \mathbf U_2(\F_p) =\F_p.\]
\begin{defi}
Assume there exists a complete flag of $(\F_p,G)$-modules \[ \nabla_3:0 \subset  V_{1} \subset V_{2} \subset V_{3},\] such that the truncated extension  of $(\F_p,G)$-modules \[ 0 \to  V_{1} \to  V_{2} \to  V_{2}/ V_{1} \to 0\] is isomorphic to $\mathcal E_x$, and such that the quotient extension  \[ 0 \to  V_{2}/ V_{1} \to  V_{3} / V_{1} \to  V_{3}/ V_{2} \to 0\] is isomorphic to $\mathcal E_y$. We then say that $\mathcal E_x$ and $\mathcal E_y$ glue, to  the complete flag $ \nabla_3$.
\end{defi}
\ni A fundamental fact is that the extensions $\mathcal E_x$ and $\mathcal E_y$ glue to a $ \nabla_3$ as above, if and only if the cup-product $$a=(x)_p \cup (y)_p \in H^2(F,\mu_p^{\otimes 2})=  H^2(F,\mu_p)=\Br(F)[p]$$ vanishes.\\
Assume that this is the case, and let $ \nabla_3$ be such a gluing. Using the same construction as above, for  coefficients in $\Z/p^2$, we get the following. \\The complete flag $ \nabla_3$ lifts to a complete flag of $(\Z/p^2,G)$-modules   \[ \nabla_{3,2}:0 \subset  V_{1,2} \subset V_{2,2} \subset V_{3,2},\] with trivial graded pieces \[L_{i,2}=\Z/p^2, \] if and only if $(x)_p$ and  $(y)_p$, respectively, lift to classes $$(X)_{p^2},(Y)_{p^2} \in  H^1(F,\mu_{p^2}),$$  such that $$(X)_{p^2} \cup (Y)_{p^2}=0 \in H^2(F,\mu_{p^2}^{\otimes 2})=  H^2(F,\mu_{p^2})=\Br(K)[p^2].$$
We now show that \begin{center}\textit{ $F$, $x$ and $y$ can be chosen, so that the liftability property above fails.} \end{center} Equivalently:
\begin{itemize}
\item{The extensions $\mathcal E_x$ and $\mathcal E_y$ glue, to  a $ \nabla_3$ as above.}
\item{The  flag $ \nabla_3$ does not admit a lift to a flag of $(G,\Z/p^2)$-modules $\nabla_{3,2}$, with trivial graded pieces.}
\end{itemize}

\ni It follows that  $ [\nabla_3] \in H^1(\Gal(F_s/F), \mathbf U_3(\F_p))$ cannot be lifted via \[H^1(\Gal(F_s/F), \mathbf U_3(\Z/p^2)) \to H^1(\Gal(F_s/F), \mathbf U_3(\F_p)),\] completing the goal of this section. The following elementary result is the key. Alternatively, one may  use a deeper, yet more involved result: \cite{Ka}, Theorem 2.1.
\begin{prop}[{\cite[Th\'eor\`eme 1]{T} or \cite[Exercise 10.5]{TW}}]\hfill  \\ 
 Let $p$ be an odd prime. Put \[F:=\C(x_1,x_2,y),\] \[x:=(x_1^p-y)(x_2^p-y) \in F\] and \[M:=F(x^{\frac 1 p}, y^{\frac 1 p}).\]Consider the cyclic  algebra \[A:=(x)_{p^2} \cup (y)_{p^2} \in  \Br(F) . \]  It is of exponent $p$, split by $M/F$. There do not exist elements $u,v \in F$ such that \[ [A]=(u)_{p} \cup (y)_{p}+(v)_{p} \cup (x)_{p} \in \Br_p(F).\]
 \end{prop}
\ni This proposition being granted, assume that $(x)_p$ and  $(y)_p$  lift to classes $$(X)_{p^2},(Y)_{p^2} \in  H^1(F,\mu_{p^2}),$$  such that $(X)_{p^2} \cup (Y)_{p^2}=0$. Write $$(X)_{p^2}=(x)_{p^2} -p (u)_{p^2}$$ and  $$(Y)_{p^2}=(y)_{p^2} -p (v)_{p^2},$$ for $u,v \in F^\times$. Expanding the equality  $(X)_{p^2} \cup (Y)_{p^2}=0$,  we get \[ [A]=(u)_{p} \cup (y)_{p}+(x)_{p} \cup (v)_{p} \in \Br(F),\] contradicting the  above.\\

\noindent We did not exclude the possibility that $\rho_1$  lifts to a representation \[G \to \B_3(\Z/p^2),\] but I believe  it does not. This would imply the non-liftability of the versal $\B_3(\F_p)$-Galois representation (over $K=\C$). From there, one could derive  an alternate proof of Theorem \ref{thmnonlift}, for $p$ odd (see section \ref{secnonliftgen} for the definition of 'versal', and more).\\

\noindent Actually, in the literature, the   first simple triangular counter-example (i.e. for $\B_3$) was given in the note \cite{F3}. Observe that its construction heavily relies on the assumption $\mu_{p^2} \nsubseteq F$, and would not work over number fields containing $\mu_{p^2}$. \\
In the sequel, write $\Gal(F)$ for $\Gal(F_s/F)$. For the sake of  concreteness, the statement given next slightly   differs from \cite{F3}, that only deals with  fields of Laurent series, i.e. with  representations of the absolute Galois groups $\Gal(F((T)))$. It is clear however, that  the construction thereof is  `non-formal': it indeed provides a $\rho_1$ as below.
\begin{prop}(See \cite{F3}.)\\
Let $p$ be an odd prime. \\ Let $F$ be a number field, such that  $\mu_p \subset F$ but $\mu_{p^2} \nsubseteq F$. \\ \ni Pick a $1$-cocycle 
\[c:  \Gal(F)\to \F_p \simeq \mu_p^{\otimes 2},\]  that does not lift to a $1$-cocycle\[  \Gal(F)\to \mu_{p^2}^{\otimes 2}.\] [From  class field theory, it is not hard to see that such a $c$ exists.]\\
 Denote by 
\[t:  \Gal(F(T))\to \mu_p \] the  $1$-cocycle (which is here a homomorphism) corresponding, via Kummer theory, to \[(T) \in H^1(F(T),\mu_p)=F(T)^\times /F(T)^{\times p}.\]
Consider  the representation  \[\rho_1:\Gal(F(T)) \to \mathbf B_3(\F_p),\]   given by the formula
\[  \begin{pmatrix}
1 &   t& t^2+c \\ 0 & 1 & 2t \\ 0 & 0 & 1
 \end{pmatrix} . \]

\ni It does not lift to a representation \[\rho_2:\Gal(F(T)) \to \B_3(\Z/p^2).\]
    In fact, it does not even lift formally at $T$: the composite representation  \[\widehat \rho_1:\Gal(F((T))) \xrightarrow{\nat} \Gal(F(T)) \xrightarrow{\rho_1} \B_3(\F_p)\] does not lift to a representation \[\Gal(F((T))) \to \B_3(\Z/p^2).\]
\end{prop}

\section{Positive results, over all fields.}\label{sec+}
\subsection{When $p$ is odd...\\}
\ni ... two-dimensional Galois representations always lift, as  the following general statement shows. It is  Theorem 6.1 of  \cite{DCF0}, whose  main contribution is to  provide a formalism  that is  insensitive to $k$.

\begin{thm}\label{liftp2d3}Let $F$ be a field. Let $k$ be a field of characteristic $p$. Then, the arrow \[H^1(\Gal(F_s/F), \GL_2(\W_2(k))) \to H^1(\Gal(F_s/F),\GL_2(k))\] is surjective.

\end{thm}

\ni It is natural to ask, if this generalises to Witt vectors of higher length.

\begin{qu}
    Let $F$ be a field. Let $k$ be a field of characteristic $p$.\\ Is the arrow \[H^1(\Gal(F_s/F), \GL_2(\W_3(k))) \to H^1(\Gal(F_s/F),\GL_2(k))\]  surjective?
\end{qu}
\ni To my knowledge, the answer is not known. A partial result is \textit{stable} liftability of two-dimensional Galois representations, again furnished by \cite{DCF0}, Theorem 6.1. Precisely, given a two-dimensional (Galois) representation $V$ over $k$, there exists another (explicit)  representation $V'$, such that $V \bigoplus V'$ lifts to a representation over $\W(k)$ (the ring of full $p$-typical Witt vectors.)\\
Proving/disproving actual liftability  to $\W_3(k)$, would definitely require some new insight.
\subsection{When $p=2$...\\}
...Galois representations over $\F_2$ lift up to dimension four, as  follows.
 \begin{prop}\label{liftp2d3}[dimension three, over $\F_2$.]\\ Let $F$ be a field of characteristic not $2$, containing $\mu_4$. Then, the arrow \[H^1(\Gal(F_s/F), \mathbf U_3(\Z/4)) \to H^1(\Gal(F_s/F), \mathbf U_3(\Z /2))\] is surjective.\\

\end{prop}
\begin{dem}
It  is a straighforward adaptation of the construction of section \ref{Heisenbergpodd}, to the case  $p=2$. This is done using the following fact, well-known to specialists (see  Proposition 5.2 of \cite{LLT}). Let $M:=F(x^{ \frac 1 2}, y^{ \frac 1 2})/F$ be a biquadratic extension. Then, every central simple algebra \[[A] \in \Br_2(M/F)\] is of the shape  \[ [A]=(u)_{2} \cup (y)_{2}+(x)_{2} \cup (v)_{2} ,\] for some $u,v \in F^\times$. Details are left to the interested reader.
\end{dem}
 \ni To my knowledge, it is unknown  whether the surjectivity statement  of  Proposition \ref{liftp2d3} holds for $\U_4$ (I believe it does not). However, there is the following result.
 \begin{prop}\;[dimension four over $\F_2$, See \cite{DCF0}, Theorem 6.1.]\label{liftp2d4} \\ Let $F$ be a field. Then, the arrows \[H^1(\Gal(F_s/F), \mathbf B_4(\Z/4)) \to H^1(\Gal(F_s/F), \mathbf B_4(\Z /2))\]  and \[H^1(\Gal(F_s/F), \GL_4(\Z/4)) \to H^1(\Gal(F_s/F), \GL_4(\Z /2))\] are surjective.
\end{prop}

\begin{rem}
 Let $k$ be a field of characteristic $2$. If $k \neq \F_2$,  Proposition \ref{liftp2d3} actually fails over $k$; see Theorem \ref{thmnonlift}.
\end{rem}

\begin{rem}
 In \cite{DCF0}, Theorem 6.1 is stated for $\GL_4$ only, but the proof clearly gives the same result for $\B_4$.
\end{rem}

\section{Negative generic results.}\label{secnonliftgen}

\ni In their recent work \cite{MeSc1} and \cite{MeSc2},  Merkurjev and Scavia gave the list of all pairs $(p,d)$, such that  Question \ref{qugallift} has an affirmative answer \textit{for all fields $F$}. They also solve the same problem, replacing   $\GL_d$ by $\B_d$, and $\F_p$ (resp. $\Z/p^2$) by an arbitrary field $k$ of characteristic $p$ (resp. $\W_2(k)$). Their result can be formulated in a simple way:  for $\GL_d$ or $\B_d$, there are no other cases in which \textit{all} Galois representations  lift mod $p^2$, than those listed in Section \ref{sec+}. Precisely:

\begin{thm}\label{thmnonlift}Let $k$ be a field of characteristic $p$, and let $d\geq 1$ be an integer.  The following are equivalent. \begin{enumerate}
    \item{ The arrow \[H^1(\Gal(F_s/F), \GL_d(\W_2(k))) \to H^1(\Gal(F_s/F), \GL_d(k))\] is surjective for all fields $F$.} \item{ If $\vert k \vert \geq 3$, then  $d \leq 2$.\\ If $\vert k \vert =2$ (equivalently, if $p=2$ and $k=\F_2$), then   $d \leq 4$.}
\end{enumerate}

\end{thm}

\ni In fact, Theorem 1.1 of  \cite{MeSc2} is slightly more precise. Here are some details. It is not too  hard to reduce to the case of a finite field $k=\F_q$. Starting with any field $K$, one can then consider the generic $d$-dimensional Galois representation over $K$, with coefficients in $k$. It consists of a field extension $F/K$, together with a representation \[\rho_{d,\vers}: \Gal(F_s/F) \to \GL_d(k)\] that is \textit{versal}. Roughly speaking, this means that, for every field exension $L/K$, and for every Galois representation \[\rho: \Gal(L_s/L) \to \GL_d(k),\]  $\rho_{d,\vers}$ has a specialisation that is conjugate to $\rho$. Equivalently, it is a \textit{versal torsor} over $K$, for the finite group $\GL_d(k)$,   in the sense of \cite{GMS}.  Liftability of all $d$-dimensional representations, over all field extensions of $K$, is thus equivalent to that of $\rho_{d,\vers}$. Therefore, one focuses on  disproving  liftability of $\rho_{d,\vers}$. A meaningful observation, is that \begin{center} \textit{Non-liftability of  $\rho_{d,\vers}$ implies that of $\rho_{d',\vers}$, for all $d'>d$.} \end{center} This is a consequence of item (1) of  Lemma \ref{lemtech}, applied to $V:=\rho_{d,\vers}$ and to the trivial representation $V'=k^{d'-d}$. At the light of this Lemma, the work is thus to disprove liftability of $\rho_{d,\vers}$, in the following cases: \begin{enumerate}
    \item{$p\geq 3$, $k=\F_{p}$, and $d =3$,} \item{$p=2$, $k=\F_{2^r}$ for $r \geq 2$,  and $d = 3$,} \item{$p=2$, $k=\F_2$ and $d =5$.}
\end{enumerate}
Item (1) is treated in \cite{MeSc1}, as an application of a result of independent interest: the  computation of \textit{negligible} cohomology classes in $H^2(G,M)$,  for all finite groups $G$ and  all finite $G$-modules $M$, over fields  containing enough roots of unity.  This result  is explicit, and its proof is nicely constructive.\\ \ni In item (2), the crucial case is $k=\F_4$. As far as I can see, this definitely requires more subtlety, than in the arguments of section \ref{Heisenbergpodd}.\\ \ni In the current version of  \cite{MeSc2}, item (3) is the hardest one- dealt with by intricate (though elementary) computations.  These can hopefully be simplified.

\section{Positive results, over specific fields.}\label{secspeclift}
\ni In \cite{KL}, Khare and Larsen prove lifting statements for   Heisenberg-Galois representations (=representations of absolute Galois groups, with values in $\U_3$), when the field $F$ is a global field, or a non-archimedean local field, containing $\mu_{p^2}$.\\ In particular, they prove the following result. At the time this survey is written, it is the only general result available in the literature, about Question \ref{qugallift} for global fields, in dimension $d \geq 3$.

\begin{prop}\;[See \cite{KL}, Theorem 5.4.]\label{liftnumber} Suppose that $p$ is odd. Let $F$ be a local field, or a number field, containing $\mu_{p^2}$.  Then, the following arrows are surjective: \[H^1(\Gal(F_s/F), \mathbf U_3(\Z/p^2)) \to H^1(\Gal(F_s/F), \mathbf U_3(\F_p))\]  and \[H^1(\Gal(F_s/F), \GL_3(\Z/p^2)) \to H^1(\Gal(F_s/F), \GL_3(\F_p)).\] 

\end{prop}
\begin{rem}
    Here again, the case of $\GL_3$ follows from that of $\U_3$, for purely group-theoretic reasons- see section \ref{remrescor}.
\end{rem}

\begin{rem}
 \ni Proposition \ref{liftnumber} still holds upon replacing $\F_p$  by a field $k$ of characteristic $p$ (and accordingly, replacing $\Z/p^2$ by $\W_2(k)$). This upgrade is at the cost of  minor modifications in proofs, as the interested reader may check. 
 
 \end{rem}
\ni For local fields, the proof of Proposition \ref{liftnumber}   just uses the following properties of $F$. \begin{enumerate}
    \item {It contains $p^2$-th roots of unity.}  \item {The $\F_p$-vector space $H^2(F,\F_p)$ is one-dimensional.}  \item {The cup-product  $$H^1(F,\F_p) \times H^1(F,\F_p) \to H^2(F,\F_p) $$ is  a perfect pairing of finite-dimensional $\F_p$-vector spaces.}
\end{enumerate}

\ni In the recent work \cite{CDF},  very general lifting theorems are proved for the so-called \textit{$p$-manageable} profinite groups. We refer to  \cite{CDF} for details on the  material presented next. 
\begin{defi}\label{Defimana}

Let $G_p$ be a pro-$p$-group. Say that $G_p$ is  $p$-manageable if the following conditions are satisfied. \begin{enumerate}

     \item {The $\F_p$-vector space $H^2(G_p,\F_p)$ is one-dimensional.}  \item {The cup-product  pairing  of (possibly infinite-dimensional) $\F_p$-vector spaces \[H^1(G_p,\F_p) \times H^1(G_p,\F_p) \to H^2(G_p,\F_p) \simeq \F_p, \]has trivial (left) kernel.}
\item{ There exists a continuous character  \[\theta_p: G_p \to \Z_p^\times,\] with the following property. Set $\Z_p(1):=\Z_p$, on which $G_p$ acts via $\theta_p$. Then, for every $r \geq 2$, the natural arrow \[H^1(G_p, (\Z/p^r)(1)) \to H^1(G_p, \F_p(1))\]is onto.\\}
\ni Let $G$ be any profinite group. Let $G_p \subset G$ be a pro-$p$-Sylow. Say that $G$ is $p$-manageable, if  $G_p$ is $p$-manageable, and the character $\theta_p$ of item (3) extends to a character \[\theta: G \to \Z_p^\times. \]
\end{enumerate}

\end{defi}

\begin{rem}
  Let $G_p$ be a $p$-manageable pro-$p$-group. Using item (2), one can prove  that the character  $\theta_p: G_p \to \Z_p^\times$ in item (3) is unique.
\end{rem}
\ni Let us give three famous examples of $p$-manageable profinite groups. 
\begin{ex}(Absolute Galois groups of local fields.)\\
Let $F$ be a finite extension of $\Q_l$, or of $\F_l((T))$, with $l=p$ allowed. Then $G:=\Gal(F_s/F)$ is $p$-manageable. Moreover, in (3), $\Z_p(1)$ is the Tate module (of roots of unity of $p$-primary order) if $\carac(F) \neq p$, or   $\Z_p(1)=\Z_p$ if $\carac(F)= p$. 
\end{ex}

\begin{ex}(Fundamental groups of curves.)\\
Let $G$ be the algebraic fundamental group of a smooth proper complex curve of genus $g>0$. Then $G$ is $p$-manageable. Moreover, in (3),  $\theta$ is trivial, i.e. $\Z_p(1)=\Z_p$.
\end{ex}

\begin{ex}(Demushkin groups.)\\
  Let $G$ be a pro-$p$-group. If $G$ satisfies items (1) and (2) of Definition \ref{Defimana}, say that $G$ is a Demushkin group. [In the classical terminology, one also requires that $G$ be finitely generated, or equivalently that $H^1(G_p,\F_p)$ is finite.]  If $G$ is finitely generated, then a character $\theta$ as in item (3) exists, so that  $G$ is $p$-manageable. For a proof that does not use dualizing modules, see    \cite{CDF}, Proposition 5.1.
\end{ex}

\ni \textit{Until the end of this section, $G$ is a $p$-manageable profinite group, and $\theta,\Z_p(1)$ are as in item (3) of Definition \ref{Defimana}.\\} 

\ni A variant of Question \ref{QuGen} for $G$, then has a very strong positive answer. \\To state it, we introduce the following notation.

\begin{defi}\label{notaV}
  For $r\geq 2$, consider a triangular representation  \[\rho_r\colon G\to\B_d(\Z/p^r).\] Denote its mod $p$ reduction  by \[\rho_1\colon G\to\B_d(\F_p).\] Denote by  $V_r$ the representation of $G$ on the free $\Z/p^r$-module $(\Z/p^r)^d$, furnished by $\rho_r$. Observe that $V_r$ is naturally equipped with a complete flag of  representations of $G$ over $\Z/p^r$. 
\end{defi}

\begin{defi}(Wound Kummer flag, for $r=1$.) \label{defiwound1}\\
     Consider a triangular representation  \[\rho_1\colon G\to\B_d(\F_p).\] Observe that $\rho_1$ gives rise to homomorphisms, for  $i=1,2,\dots d-1$, \[\lambda_i:G\to\B_2(\F_p),\]  corresponding matrixwise, to the  $(2 \times 2)$ blocks centered at the diagonal. \\If one of the following three equivalent conditions is satisfied, we say  that the (triangular) representation $\rho_1$ is wound Kummer. \begin{enumerate}
         \item{For all $i$, $p$ divides $\vert \Im(\lambda_i) \vert$.}  \item{For all $i$, the  extension (of one-dimensional representations of $G$ over $\F_p$) corresponding to  $\lambda_i$ is  non-split.}  \item{There is a unique $G$-invariant complete flag on $V_1$, given by $\rho_1$.}
     \end{enumerate} 
\end{defi}
\begin{defi}(Wound Kummer flag, general case.) \label{defiwound}\\
   For $r\geq 1$,  consider a triangular representation  \[\rho_r\colon G\to\B_d(\Z/p^r).\]  The diagonal of $\rho_r$ gives $d$ multiplicative characters (for  $i=1,2,\dots d$) \[\chi_i:G\to (\Z/p^r)^\times.\]  Say  that $\rho_r$ is wound Kummer if the following conditions hold. \begin{enumerate}
         \item{ $\rho_1$ is wound Kummer, in the sense of Definition \ref{defiwound1}.} \item{For $i=1,\ldots,d$,  the order of the character \[\chi_i.\theta^i:G\to (\Z/p^r)^\times \]  divides $(p-1)$.}
     \end{enumerate} 
\end{defi}

\begin{rem}
 Item (2) of Definition \ref{defiwound} can be reformulated as: the $p$-primary parts of $\chi_i$ and (the reduction modulo $p^r$ of) $\theta^{-i}$ coincide.
\end{rem}

\ni We do not provide here a precise definition of a Kummer representation, but the rough idea goes like this.   A representation \[\rho_r\colon G\to\U_d(\Z/p^r)\]  is Kummer, if the combinatorics of partial splittings of $\rho_r$, is the same as that of its mod $p$ reduction $\rho_1$.\\

\ni A   very strong step-by-step lifting result is available for wound Kummer, resp. for Kummer representations:  Theorems 7.10, resp. 7.21 of \cite{CDF}. In this survey, we contend ourselves with two corollaries of these: Theorems \ref{thmKum} and  \ref{prop12} below. Separately, they illustrate  the two meanings of ``step-by-step'': \begin{itemize}
    \item{ On the one hand,   w.r.t. the depth $r$ of the flag- see Theorem \ref{thmKum}.} \item{On the other hand,  w.r.t. the  dimension $d$- see item (2) of  Theorem \ref{prop12}.}
\end{itemize}  
\ni The construction of the liftings is explicit,  using an iterative deformation process, along which extensions are manipulated via elementary operations: Baer sum, push-forward and pull-back. A key technique is a mixture of gluing and lifting, called ``gluifting'' (\cite{CDF}, section 5). 
\begin{thm}\label{thmKum}
Let $d,r \geq 1$ be  integers. Let $\rho_r$ be a $d$-dimensional wound Kummer, resp. a Kummer representation. In the Kummer case, assume that $\Z/p^{r+1}(1)=\Z/p^{r+1}$. Then $\rho_r$  lifts, to a wound Kummer, resp. to a Kummer representation $\rho_{r+1}$.
\end{thm}

\begin{rem}
Kummer representations provide a natural framework, where  one can lift certain Galois representations via the arrow  \[ H^1(\Gal(F_s/F),\GL_d(\Z/p^{r+1})) \to  H^1(\Gal(F_s/F),\GL_d(\Z/p^r)), \]  for $d \geq 1$ and $r \geq 2$. Results of this kind are very rare. Observe that, for general $F$, this arrow is never surjective. Here is an example. Take $r=2$ if $p$ is odd, or $r=3$ if $p=2$. Then, surjectivity   fails already for $d=1$. Indeed, it is not hard to check, that it is equivalent to surjectivity of \[ H^1(\Gal(F_s/F),\Z/p^2) \to  H^1(\Gal(F_s/F),\F_p), \] which does not hold in general, e.g. for $F=\mathbf Q$. For fields containing $\mathbb C$, there should be counter-examples also for $d=2$, but I do not have any in mind. 
\end{rem}

\begin{thm}\;[\cite{CDF}, Corollary 7.22.]\label{prop12}\\
Consider a representation \[\rho_1 \colon G\to\U_d(\F_p). \] Let $r \geq 2$, and assume that  $\Z/p^r(1)=\Z/p^r$. Then, the following hold. \begin{enumerate}
    \item{There exists a lift  of $\rho_1$, to a  representation   \[\rho_r\colon G\to\U_d(\Z/p^r).\] } \item{Futhermore, $\rho_r$ can be picked such that the natural map \[ H^1(G,V_r) \to H^1(G,V_1)\] is surjective. }
\end{enumerate}
\end{thm}

\begin{rem}
   For  item (1) to hold, the condition $\Z/p^r(1)=\Z/p^r$ is necessary.\end{rem}
    
    \begin{rem}
For absolute Galois groups of local fields, it follows from results of \cite{EG}, that item (1) holds  with $\B_d$ in place of $\U_d$,  without the assumption $\Z/p^r(1)=\Z/p^r$. However, by \cite{CDF}, Example 7.23, this cannot be achieved using powers of the cyclotomic character for the diagonal of $\rho_r$. In other terms, the $d$ one-dimensional graded pieces of $V_r$ cannot in general be  of the  form $L_{i,r}:=\Z/p^r(n_i)$, for $i=1,\ldots d$.
\end{rem}
 \begin{rem}
      It is not clear to me how to prove item (2) of Theorem \ref{prop12}, using the material of \cite{EG}.  
 \end{rem}
\begin{rem}\label{remKumd}
Item (2) can be thought of as  ``higher-dimensional Kummer theory'', where  the $d$-dimensional representation $V_r$ replaces the one-dimensional cyclotomic module $\Z/p^r(1)$. Meanwhile,  this analogy is seriously limited: in general, there is no choice of $\rho_r$,  such that  \[ H^1(H,V_r) \to H^1(H,V_1)\] is surjective \textit{for all open subgroups} $H \subset G$.\\ Details are provided the last section of this paper.
\end{rem}

\section{There is no ``naive'' two-dimensional Kummer theory.}\label{secnonKum}

\ni In the literature, the next result is new.
\begin{prop}\label{niceexo}
 Let $F$ be a field of characteristic not $p$, containing all $p^2$-th roots of unity. Set $G:=\Gal(F_s/F)$. Assume that, for every open subgroup $H \subset G$, the cup-product pairing \[H^1(H,\F_p) \times H^1(H,\F_p) \to H^2(H,\F_p)\] is non-degenerate, meaning that its left kernel is trivial. [This is the case if $F$ is a local field.  By Lemma 3.6 of \cite{CDF},  it is also the case if $F$ is infinite and finitely generated over its prime subfield.] \\ Consider a two-dimensional representation  \[\rho_1 \colon G\to\GL_2(\F_p), \] such that $\vert \Im(\rho_1) \vert$ is divisible by $p$. \\ Then, there  does not exist a lift of $\rho_1$ to a representation  \[\rho_2 \colon G\to\GL_2(\Z/p^2), \] such that, with notation of Definition \ref{notaV}, the natural map \[H^1(H,V_2) \to H^1(H,V_1)\] is surjective for every open subgroup $H \subset G$. 

\end{prop}

\begin{dem}
 For the sake of contradiction,   assume that such a $\rho_2$ exists. By assumption, there exists an open subgroup $G' \subset G$, such that $\rho_1(G')$ is cyclic of order $p$.  Replacing $G$ by $G'$, and replacing $\rho$ by a conjugate representation, one thus reduces to the case where $\rho_1 \neq 1$ takes values onto $\U_2(\F_p)=\Z/p$. Set $H:=\Ker(\rho_1)$. If $\rho_2(H)=1$, then  the matrix \[ \begin{pmatrix}
1 & 1\\
  0& 1
\end{pmatrix} \in \GL_2(\F_p)\]
would lift (via $\rho_2$) to an element of  $\GL_2(\Z/p^2)$ of order $p$. This is absurd.  Consequently, it suffices to prove that $\rho_2(H)=1$. Since $V_2$ is a free $\Z/p^2$-module, the  natural  exact sequence of $H$-modules  \[ 0 \to pV_2 \to V_2 \to V_2/pV_2  \to 0\] reads as \[ 0 \to V_1 \to V_2 \to V_1\to 0,\] i.e. \[ (E): 0 \to \F_p^2 \to V_2  \xrightarrow{q}\F_p^2 \to 0.\] By assumption, the arrow \[q_*: H^1(H,V_2) \to H^1(H,\F_p^2)\] is onto. Since $\mu_{p^2} \subset F$, Kummer theory implies that the arrow \[H^1(H,\Z/p^2) \to H^1(H,\F_p)\] is onto. Define \[V'_2:=(\Z/p^2)^2.\] It is another (obvious) lift of the trivial representation $H$ over $\F_p$, to a representation of $H$ over $\Z/p^2$. There is  a natural extension of (trivial)  $H$-modules  \[ (E'): 0 \to \F_p^2 \to V'_2  \xrightarrow{q'}\F_p^2 \to 0,\] and by Kummer theory recalled above,  \[q'_*: H^1(H,V'_2) \to H^1(H,\F_p^2)\] is onto as well. Consider $(E)$ and $(E')$ as extensions of $\Z/p^2$-modules (with an action of $H$). As such, form their Baer difference \[ (E)-(E'):  0 \to \F_p^2 \to D \xrightarrow{\pi} \F_p^2 \to 0,\] and denote it by $\Delta$. Since $q_*$ and $q'_*$ are onto, so is the arrow \[\pi_*: H^1(H,D) \to H^1(H,\F_p^2).\] To check this, observe first that surjectivity of $q_*$ amounts to vanishing  of the connecting map associated to the extension $(E)$, reading as \[ \beta_{E}: H^1(H,\F_p^2) \to H^2(H,\F_p^2). \]The same fact holds for $q'_*$. Using that the formation of Baer sum is compatible to connecting maps, it follows that \[\beta_\Delta: H^1(H,\F_p^2) \to H^2(H,\F_p^2)\] vanishes, whence the sought-for surjectivity. Next, since $V_2$ and $V'_2$ are both lifts of the $\F_p$-vector space $V_1(=\F_p^2)$ to free $\Z/p^2$-modules, it follows that $D$ is in fact an  $\F_p$-vector space. The proof is elementary, and left to the reader. [One may use the connecting arrow $\kappa$, introduced in \cite{DCFL}, section 3.5.] Consequently, the extension of $(\F_p,H)$-modules $\Delta$ is provided by a matrix  \[ \begin{pmatrix}
D_{1,1} & D_{1,2} \\
  D_{2,1} &D_{2,2} 
\end{pmatrix},\] where  each $D_{i,j}$ is an extension of  $(\F_p,H)$-modules, of the shape \[D_{i,j}: 0 \to \F_p \to P_{i,j} \to \F_p \to 0. \]
Denote by \[ d_{i,j}=[D_{i,j}] \in H^1(H,\F_p)\] the cohomology class of $D_{i,j}$. The connecting map  \[\beta_\Delta: H^1(H,\F_p)^2 \to H^2(H,\F_p)^2\] is thus given (up to sign) by the formula \[(x_1,x_2) \mapsto (x_1 \cup d_{1,1}+ x_2 \cup d_{2,1} \;, \; x_1 \cup d_{1,2}+ x_2 \cup d_{2,2}). \] By the above, this connecting map identically  vanishes, which implies that the four $d_{i,j}$'s lie in the kernel of the cup-product pairing  \[H^1(H,\F_p) \times H^1(H,\F_p) \to H^2(H,\F_p). \] By the non-degeneracy assumption, the $d_{i,j}$'s  vanish. Consequently, the extension  $\Delta$ is trivial, which implies that $V_2 \simeq V'_2$, as $(\Z/p^2,H)$-modules. Equivalently, $H$ acts trivially on $V_2$, as was to be shown.
\end{dem}
\section{Acknowledgements.}
Thanks to Charles De Clercq and Federico Scavia for their reading and remarks.
\bibliographystyle{plain}
\bibliography{main.bib}
\end{document}